\documentclass{amsart}

\usepackage{amsmath,amsthm,amsfonts,latexsym}
\numberwithin{equation}{section}

\def\R{\mathbb{R}}
\def\Q{\mathbf{Q}}

\def\N{\mathbb{N}}
\def\K{\mathbf{K}}

\def\L{\mathbf{L}}

\def\I{\mathbf{I}}

\def\Six{\mathbf{\Sigma^2_x}}
\def\Sixy{\mathbf{\Sigma^2_{xy}}}
\def\Siyz{\mathbf{\Sigma^2_{yz}}}

\newtheorem{thm}{Theorem}[section]

\theoremstyle{definition}

\theoremstyle{remark}

\title{A Positivstellensatz which preserves the 
coupling pattern of variables}

\author{Jean B. Lasserre}

\address{LAAS-CNRS and Institute of Mathematics, ; 90C25
LAAS, 7 Avenue du Colonel Roche, 31077 Toulouse C\'edex 4, France}

\email{lasserre@laas.fr}

\date{}

\begin{document}

\begin{abstract}
We specialize Schm\"udgen's Positivstellensatz and its
Putinar and Jacobi and Prestel  refinement, to the case
of a polynomial $f\in\R[X,Y]+\R[Y,Z]$, positive on a compact basic semi-algebraic set
$\K$ described by polynomials in $\R[X,Y]$ and $\R[Y,Z]$ only, or in
$\R[X]$ and $\R[Y,Z]$ only (i.e. $\K$ is cartesian product).
In particular, we show that the preordering $P(g,h)$ (resp. quadratic module $Q(g,h)$) generated
by the polynomials $\{g_j\}\subset\R[X,Y]$ 
and $\{h_k\}\subset\R[Y,Z]$ that describe $\K$,
is replaced with $P(g)+P(h)$ (resp. $Q(g)+Q(h)$), so that the absence of coupling
between $X$ and $Z$ is also preserved in the representation.
A similar result applies with Krivine's Positivstellensatz involving the cone
generated by $\{g_j,h_k\}$. 
\end{abstract}

\maketitle




\section{Introduction}

Schm\"udgen's Positivstellensatz  \cite{schmudgen} and its 
Putinar \cite{putinar} and Jacobi and Prestel \cite{jacobi} refinements, are
very useful tools with (relatively recent) particular impact in polynomial 
optimization; see e.g. Lasserre \cite{siopt1}, Schweighofer \cite{markus}.

An interesting issue, and not only from a computational viewpoint, is to
derive a Positivstellensatz that preserves a possible coupling pattern of variables
present in the original polynomial $f$ (positive on a compact semi-algebraic set).
By this we mean that if there is no coupling of variables
$X_i$ and $X_j$ in $f$ as well as in the polynomials that describe the
basic semi-algebraic set, we would like  to obtain a representation in which the same 
property holds.

A first positive result in this vein is derived in Lasserre \cite{largescale} under some condition of the coupling pattern, known as the {\it running intersection property} in graph theory.
Such specialized representations are particularly important from a computational viewpoint, as evidenced by the impressive computational experiments presented in Waki et al. \cite{waki}, when used in polynomial optimization problems with structured sparsity.

Here we present two specialized Positivstellensatz
when $f\in\R[X,Y]+\R[Y,Z]$, that is,
when there is {\it no} coupling between variables $X$ and $Z$ in $f$.
While in the first one, the compact basic semi-algebraic set $\K$ is also described
by polynomials in  $R[X,Y]$ and $\R[Y,Z]$ only, in the second one 
$\K$ is described by polynomials in $\R[X]$ and $R[Y,Z]$ only, i.e., $\K=\K_x\times\K_{yz}$ for some $\K_x,\K_{yz}$.
Our result  does {\it not} require any assumption at all  (except of course compactness of $\K$).
 Although part of our result in the first case, namely Theorem \ref{thmain}(b) below, could be derived from  \cite{largescale} with appropriate modifications, 
 the general form of our specialized Positivstellensatz in Theorem \ref{thmain}(a)-(b)-(c) is not
 apparent from \cite{largescale}, and we think it is 
important enough and of self-interest to deserve a special treatment and presentation
to an audience not necessarily aware of the more computational oriented result \cite{largescale}. On the other hand, the more involved case where $\K$ is a cartesian product cannot be deduced from \cite{largescale}.
\vspace{0.2cm}

\noindent
{\bf Contribution.}
Let $\K_{xy}\subset\R^{n+m}$, $\K_{yz}\subset\R^{m+p}$, and 
$\K\subset\R^{n+m+p}$ be basic compact semi-algebraic sets defined by
\begin{eqnarray}
\label{kxky}
\K_{xy}&=&\{\:(x,y)\,\in\R^{n+m}\:\vert\quad g_j(x,y)\geq0,\quad j\in \I_{xy}\}\\
\label{kyz}
\K_{yz}&=&\{\:(y,z)\,\in\R^{m+p}\:\vert\quad h_k(y,z)\geq0,\quad k\in \I_{yz}\}\\
\label{setk}
\K&=&\{\:(x,y,z)\,\in\R^{n+m+p}\:\vert\quad (x,y)\in\K_{xy};\quad
(y,z)\in\K_{yz}\:\}
\end{eqnarray}
for some polynomials $\{g_j\}\subset\R[X,Y]$, $\{h_k\}\subset\R[Y,Z]$, and some finite index sets $\I_{xy},\I_{yz}\subset\N$.

Let $P(g)\subset\R[X,Y]$ and
$P(h)\subset\R[Y,Z]$ be the preordering generated by $\{g_j\}_{j\in\I_{xy}}$
and $\{h_k\}_{k\in\I_{yz}}$, respectively.  

Let $f\in\R[X,Y]+\R[Y,Z]$, i.e., there is no coupling of variables $X$ and $Z$ in $f$.

$\bullet$ We first obtain the following specialized Positivstellensatz.
 \begin{equation}
\label{res1}
\left[\,f\in\,\R[X,Y]+\R[Y,Z]\mbox{ and }f\,>\,0\mbox{ on }\K\,\right]
\quad\Rightarrow\quad f\in \:P(g)+P(h),
\end{equation}
 to compare with Schm\"udgen's Positivstellensatz which states that
$f\in P(g,h)$. 

The Positivstellensatz  (\ref{res1}) is a specialization of Schm\"udgen's Positivstellensatz where the preordering $P(g,h)$ is replaced with $P(g)+P(h)$.
And so, only polynomials in $\R[X,Y]$ and $\R[Y,Z]$ are involved in the representation (\ref{res1}). In other words, the absence of coupling between the variables $X$ and $Z$ is preserved in the Positivstellensatz.

If in addition to be compact, $\K_{xy}$ and/or $\K_{yz}$ satisfy 
Putinar's condition in \cite{putinar}, then in
(\ref{res1}) one may replace $P(g)$ and/or  $P(h)$ with the quadratic modules
$Q(g)\subset\R[X,Y]$ and $Q(h)\subset\R[Y,Z]$ generated by $\{g_j\}_{j\in\I_{xy}}$
and $\{h_k\}_{k\in\I_{yz}}$, respectively.\\

Finally, assume that $0\leq g_j\leq 1$  on $\K_{xy}$
for all $j\in \I_{xy}$, and $0\leq h_k\leq 1$ on $\K_{yz}$ for all $k\in\I_{yz}$).
If the families $\{0,1,\{g_j\}\}$ and $\{0,1,\{h_k\}\}$ generate the 
algebra $\R[X,Y]$ and $\R[Y,Z]$ respectively, then one also obtains the alternative representation
\begin{equation}
\label{res2}
\left[\,f\in\,\R[X,Y]+\R[Y,Z]\mbox{ and }f\,>\,0\mbox{ on }\K\,\right]\,
\Rightarrow\, f\in \:C(g,1-g)+C(h,1-h),
\end{equation}
where $C(g,1-g)$ (resp. $C(h,1-h)$) is the {\it cone} generated by the 
polynomials $\{g_j,1-g_j\}$ (resp. $\{h_k,1-h_k\}$).

The Positivstellensatz  (\ref{res2}) is a specialization of 
Krivine \cite{krivine} and Vasilescu
\cite{vasilescu} Positivstellensatz where the cone $C(g,h,1-g,1-h)$ is replaced with $C(g,1-g)+
C(h,1-h)$. And so, only polynomials in $\R[X,Y]$ and $\R[Y,Z]$ are involved in the representation (\ref{res2}). In other words, the absence of coupling between the variables $X$ and $Z$ is also preserved in the Positivstellensatz (\ref{res2}).

$\bullet$ When $\K_{xy}$ is now replaced with $\K_x\subset\R^n$, that is, $\K=\K_x\times\K_{yz}$, then with $f\in\R[X,Y]+\R[Y,Z]$, we now obtain
\begin{equation}
\label{res1bis}
\left[f\,>\,0\mbox{ on }\K_x\times\K_{yz}\,\right]
\quad\Rightarrow\quad f\in \:\Sixy +P(g)+P(h),
\end{equation}
where now $P(g)\subset\R[X]$ and $\Sixy\subset\R[X,Y]$ denotes the set of sums of squares. This case is more involved and cannot be derived from \cite{largescale}.

The paper is organized as follows. Our two results are stated in the next section, and
for clarity of exposition, their proofs are postponed to section \S \ref{proofs}.

\section{Main result}

Let $\R[X,Y,Z]$ denote the ring of real polynomial in the variables
$(X_1,\ldots,X_n)$,
$(Y_1,\ldots,Y_m)$ and $(Z_1,\ldots,Z_p)$. Let $\Vert x\Vert$ denote the
euclidean norm of $x\in\R^n$.

Let $\Sixy\subset\R[X,Y]$ (resp. $\Siyz\subset\R[Y,Z]$) be the space 
of elements of $\R[X,Y]$ 
(resp. $\R[Y,Z]$) that are sum of squares (in short s.o.s.).

Given a family $\{g_j\}_{j\in\I_{xy}}\in \R[X,Y]$ (resp. $\{h_k\}_{k\in\I_{yz}}\subset\R[Y,Z]$) for some finite index set $\I_{xy}\subset\N$ (resp. $\I_{yz}\subset\N$), denote by $P(g)$ (resp. $P(h)$) the {\it preordering} 
generated by $\{g_j\}$ (resp. $\{h_k\}$). That is, $\sigma\in P(g)$ if
\begin{equation}
\label{preordering}
\sigma\,=\,\sum_{J\subseteq \I_{xy}}\,\sigma_J\,g_J\,(=\,
\sum_{J\subseteq \I_{xy}}\,\sigma_J\prod_{j\in J}g_j)
\quad\mbox{with}\quad \sigma_J\in\Sixy\quad\forall\,J\subseteq\I_{xy}
\end{equation}
(with the convention $g_J:=\prod_{j\in J}g_j\equiv 1$ if $J=\emptyset$), 
and same thing for $P(h)$ with obvious adjustments.
Similarly, denote by $Q(g)\subset\R[X,Y]$ the {\it quadratic module} generated by $\{g_j\}_{j\in \I_{xy}}\subset\R[X,Y]$. That is, $\sigma\in Q(g)$ if
\begin{equation}
\label{module}
\sigma\,=\,\sigma_0+\sum_{j\in  \I_x}\,\sigma_j\,g_j\quad\mbox{with}\quad 
\sigma_j\in\Sixy\quad\forall\,j\in\I_{xy}\cup\{0\},
\end{equation}
and same thing for $P(h)$ with obvious adjustments.\\

Finally, denote by $C(g,1-g)$ (resp. $C(h,1-h)$) the {\it cone} generated by the family
$\{g_j,1-g_j\}$ (resp.$\{h_k,1-h_k\}$), i.e. $\sigma\in C(g,1-g)$ if
\begin{equation}
\label{cone}
\sigma\,=\,\sum_{\alpha,\beta\in\N^{\vert \I_{xy}\vert}}\,
c_{\alpha\beta}\,g^{\alpha}\,(1-g)^{\beta}\,=\,
\sum_{\alpha,\beta\in\N^{\vert \I_{xy}\vert}}\,
c_{\alpha\beta}\prod_{j\in\I_{xy}}\,g_j^{\alpha_j}\,(1-g_j)^{\beta_j} 
\end{equation}
for some {\it nonnegative} scalar coefficients $\{c_{\alpha\beta}\}$.\\


\begin{thm}
\label{thmain}
Let $\K_{xy}\subset\R^n$, $\K_{yz}\subset\R^m$, and $\K\subset \R^{n+m+p}$ be the basic compact semi-algebraic sets defined in (\ref{kxky})-(\ref{setk}), and assume that
$\K$ has nonempty interior. Let
$f\in \R[X,Y]+\R[Y,Z]$.\\

{\rm (a)} If $f$ is positive on $\K$ then $f\in P(g)+P(h)$.\\

{\rm (b)} If $N-\Vert (X,Y)\Vert^2\in Q(g)$ and/or $N-\Vert (Y,Z)\Vert^2\in Q(h)$ for some scalar $N$, and if  $f$ is positive on $\K$,
then in (a) one may replace $P(g)$ with $Q(g)$ and/or $P(h)$ with $Q(h)$. \\

{\rm (c)} Assume that the 
$g_j$'s and $h_k$'s are normalized, i.e.,
$0\leq  g_j\leq 1$ on $\K_{xy}$ for all $j\in\I_{xy}$, and
$0\leq  h_k\leq 1$ on $\K_{yz}$ for all $k\in\I_{yz}$. In addition, assume that the family
$\{0,1,g_j,j\in\I_{xy}\}$ (resp. $\{0,1,h_k,k\in\I_{yz}\}$)
generates the algebra $\R[X,Y]$ (resp. $\R[Y,Z]$).
If  $f$ is positive on $\K$ then $f\in C(g,1-g)+C(h,1-h)$.
\end{thm}

The proof is postponed to \S \ref{proofthmain}. One may see that in
the various representations of $f$ in Theorem \ref{thmain}(a)-(b)-(c), all
polynomials in $P(g),P(h)$, $Q(g),Q(h)$, $C(g,1-g)$ and $C(h,1-h)$ are either in $\R[X,Y]$
or $\R[Y,Z]$, and so, there is no coupling of variables $X$ and $Z$, as in $f$. In other words, the coupling pattern of variables is preserved in each representation.

As already mentioned, when $N-\Vert (X,Y)\Vert^2\in Q(g)$ and $N-\Vert (Y,Z)\Vert^2\in Q(h)$ for some scalar $N$, then Theorem \ref{thmain}(b) can be derived from a result stated in \cite{largescale}, with appropriate modifications. This is because the condition  (1.3) in \cite{largescale}
(known as the running intersection property in graph theory)
is satisfied, by taking $I_1:=(X,Y)$, $I_2=(Y,Z)$.

We next consider the more involved case where 
$\K$ is the cartesian product $\K_x\times\K_{yz}$, with
\begin{equation}
\label{kx}
\K_x\,=\,\{\:x\in\R^n\::\quad g_j(x)\,\geq\,0,\quad j\in \I_x\:\},
\end{equation}
for some polynomials $\{g_j\}\subset\R[X]$. Now, both sets
$P(g)$ and $Q(g)$ are subsets of $\R[X]$.
\vspace{0.2cm}

\begin{thm}
\label{coro1}
Let $\K_{yz}\subset\R^{m+p}$ and $\K_{x}\subset\R^{n}$ be compact basic semi-algebraic sets as defined in
(\ref{kyz}) and (\ref{kx}) respectively. Assume that 
$\K:=\K_x\times\K_{yz}\subset \R^{n+m+p}$ has nonempty interior, and let 
$f\in \R[X,Y]+\R[Y,Z]$.\\
 
{\rm (a)}  If $f$ is positive on $\K$ then $f\in \Sixy+P(g)+P(h)$.\\

{\rm (b)} If $N-\Vert X\Vert^2\in Q(g)$ and/or $N-\Vert (Y,Z)\Vert^2\in Q(h)$ for some scalar $N$, and if  $f$ is positive on $\K$,
then in (a) one may replace $P(g)$ with $Q(g)$ and/or $P(h)$ with $Q(h)$. \\
\end{thm}

For a proof see \S \ref{proofcoro1}. So if we compare Theorem \ref{thmain}(a)-(b) and
Theorem \ref{coro1}(a)-(b) we see that in the latter the preordering 
$P(g)$ and the quadratic module $Q(g)$ are now in $\R[X]$ instead of
$\R[X,Y]$ in the former. On the other hand we need some term of $\Sixy$.
Finally, Theorem \ref{coro1}(b) cannot be deduced from \cite{largescale}.

\section{Proofs}
\label{proofs}
We first need to introduce some additonal notation and definitions.

\subsection{Notation and definitions}
\label{notation}
Let $u=(u_{\alpha\beta\gamma})$ be a sequence indexed in the canonical 
basis $(X^\alpha Y^\beta Z^\gamma)$ of $\R[X,Y,Z]$, and define the linear functional
$L_u:\R[X,Y,Z]\to\R$ to be:
\[f \:(=\sum_{\alpha,\in\N^{n}, \beta\in\N^{m},\gamma\in\N^p} f_{\alpha\beta\gamma}\,
X^\alpha\,Y^\beta\,Z^\gamma)\:\mapsto\: L_u(f) =
\sum_{\alpha,\in\N^{n}, \beta\in\N^{m},\gamma\in\N^{p}} 
f_{\alpha\beta\gamma}\,u_{\alpha\beta\gamma}.\]

\noindent
{\bf Moment matrix.} With a sequence $u=(u_{\alpha\beta\gamma})$ indexed in the canonical 
basis of $\R[X,Y,Z]$ is associated the {\it moment matrix} $M_r(u)$ of order $r$, and defined by 
\[M_r(u)((\alpha,\beta,\gamma),(\alpha',\beta',\gamma'))\,=\,
L_u(X^{\alpha+\alpha'} Y^{\beta+\beta'} Z^{\gamma+\gamma'})\,=\,
u_{\alpha+\alpha',\beta+\beta',\gamma+\gamma'},\]
for all $\alpha,\alpha'\in\N^n,\beta,\beta'\in\N^{m},
\gamma,\gamma'\in\N^p$, with $\vert\alpha+\beta+\gamma\vert\leq r$,
and $\vert\alpha'+\beta'+\gamma'\vert\leq r$.

An infinite sequence $u=(u_{\alpha\beta\gamma})$ has a 
representing measure $\mu$ if
\[u_{\alpha\beta\gamma}\,=\,\int X^\alpha Y^\beta Z^\gamma\,d\mu (X,Y,Z),
\qquad\alpha\in\N^{n},\,\beta\in\N^{m},\,\gamma\in\N^p.\]
Let $\R_r[X,Y,Z]\subset\R[X,Y,Z]$ be the space of polynomials of degree at most $r$, and denote by $s(r):={n+m+p+r\choose r}$ its dimension. If 
$\mathbf{f}=(f_{\alpha\beta\gamma})\in\R^{s(r)}$ denotes 
the vector of coefficients of an arbitrary polynomial  $f\in\R_r[X,Y,Z]$ then
\begin{equation}
\label{mom1}
\langle \mathbf{f},M_r(u)\mathbf{f}\rangle\,=\,L_u(f^2)\,=\,\int f^2\,d\mu\geq0,
\end{equation}
and so, as $f\in\R_r[X,Y,Z]$ was arbitrary, $M_r(u)\succeq0$. 

We next denote by $M_r(u,xy)$ the moment submatrix obtained from $M_r(u)$ by
retaining only those rows and columns $(\alpha,\beta,\gamma)$ with $\gamma=0$.
Similarly, denote by $M_r(u,yz)$ the moment submatrix obtained from $M_r(u)$ by
retaining only those columns and rows $(\alpha,\beta,\gamma)$ with $\alpha=0$.
Introducing the subsequences
$u^{xy}:=(u_{\alpha\beta 0})$ and $u^{yz}:=(u_{0\alpha\beta})$, notice that
$M_r(u, xy)$ is just the moment matrix $M_r(u^{xy})$ of the sequence
$u^{xy}$ indexed in the canonical basis $(X^\alpha Y^\beta)$ of $\R[X,Y]$. Similarly,
$M_r(u, yz)$ is the moment matrix $M_r(u^{yz})$ of the sequence
$u^{yz}$ indexed in the canonical basis $(Y^\alpha Z^\beta)$ of $\R[Y,Z]$. 
\vspace{0.2cm}

\noindent
{\bf Localizing matrix.} Let $g\in\R[X,Y,Z]$ and let $u=(u_{\alpha\beta\gamma})$ be a sequence indexed in the canonical basis of $\R[X,Y,Z]$. 
The {\it localizing} matrix $M_r(g\,u)$
associated with $g$ $(=\sum_{abc}g_{abc}X^aY^bZ^c)$ and $u$, is obtained from $M_r(u)$ by

\begin{eqnarray*}
M_r(g\,u)((\alpha,\beta,\gamma),(\alpha',\beta',\gamma'))&=&
L_u( g(X,Y,Z)\, X^{\alpha+\alpha'}Y^{\beta+\beta'}Z^{\gamma+\gamma'})\\
&=&\sum_{a\in\N^{n},b\in\N^{m},c\in\N^p}g_{abc}\,u_{a+\alpha+\alpha',b+\beta+\beta',c+\gamma+\gamma'},
\end{eqnarray*}
for all $\alpha,\alpha'\in\N^{n}$, $\beta,\beta'\in\N^{m}$ and 
$\gamma,\gamma'\in\N^p$ with $\vert\alpha+\beta+\gamma\vert\leq r$
and $\vert\alpha'+\beta'+\gamma'\vert\leq r$.

Again, if $u$ has a representing measure $\mu$ then
\begin{equation}
\label{mom2}
\langle \mathbf{f},M_r(g\,u)\mathbf{f}\rangle\,=\,L_u(g\,f^2)\,=\,\int g\,f^2\,d\mu,\qquad \forall\,f\in\R_r[X,Y,Z].
\end{equation}
If $\mu$ has its support in the level set
$\{(x,y,z)\in\R^{n+m+p}\:\vert\:g(x,y,z)\geq0\}$ then $M_r(g\,u)\succeq0$.

As for the moment matrix,
one may also define the localizing matrices
$M_r(g\,u,xy)$ and $M_r(h\,u,yz)$ associated with $u$ and
$g\in\R[X,Y]$, $h\in\R[Y,Z]$, respectively. They are obtained 
from $M_r(g\,u)$ (resp. $M_r(h\,u)$) by retaining only those rows and 
columns $(\alpha,\beta,\gamma)$ with
$\gamma=0$ (resp. $\alpha=0$). They can also be considered as the 
localizing matrix $M_r(g\,u^{xy})$ (resp. $M_r(h\, u^{yz})$) associated with $g$ (resp. $h$) and the subsequence $u^{xy}$
(resp. $u^{yz}$).

\subsection{Proof of Theorem \ref{thmain}}
\label{proofthmain}
(a) For every $J\subseteq \I_{xy}$, and depending on parity, let
$2r(g_J)$ or $2r(g_J)-1$ be the degree of $g_J$ (and similarly for
$2r(h_K)$ or $2r(h_K)-1$). Let $2r_0:=\max[{\rm deg}\,f,\, \displaystyle 2\max_{J,K}[r(g_J),r(h_K)]]$, and for
$r\geq r_0$, consider the following optimization problem:

\begin{equation}
\label{pbp}
\Q_r:\quad\left\{\begin{array}{ll}
\displaystyle\min_{u} &L_u(f)\\
\mbox{s.t.}&
M_{r-r(g_J)}(g_J\,u,xy)\succeq0,\quad \forall\,J\subseteq\I_{xy}\\
&M_{r-r(h_k)}(h_K\,u,yz)\succeq0,\quad \forall\,K\subseteq\I_{yz}\\
&u_0=1
\end{array}\right.,
\end{equation}
where $g_J:=\prod_{j\in J}g_j$ for every $J\subseteq \I_{xy}$ and
$h_K:=\prod_{k\in K}h_k$ for every $K\subseteq \I_{yz}$. Recall that
$g_\emptyset\equiv1$ and so, $M_r(g_\emptyset\,u,xy)$ is just the moment matrix
$M_r(u,xy)$ defined in \S \ref{notation}; the same holds true for $h_\emptyset$.
Denote by $\inf\Q_r$ the optimal value of $\Q_r$.\\

$\Q_r$ is a convex optimization problem called a {\it semidefinite programming}
problem. Up to arbitrary fixed precision, it can be solved in time polynomial in the input size of the problem data, and efficient specialized software packages are available.
For more details the interested reader is referred to e.g. Vandenberghe and Boyd \cite{boyd}.
The dual problem $\Q^*_r$ of $\Q_r$ is also a semidefinite program, which reads
\begin{equation}
\label{pbpdual}
\Q^*_r:\quad\left\{\begin{array}{ll}
\displaystyle\max_{\{\sigma_J\},\{\psi_K\},\lambda} &\lambda\\
\mbox{s.t.}&
f-\lambda\,=\,\displaystyle\sum_{J\subseteq\I_{xy}}\,\sigma_J\,g_J
+\displaystyle\sum_{K\subseteq\I_{yz}}\,\psi_K\,h_K\\
&\\
&\{\sigma_J\}\subset\Sixy,\quad\{\psi_K\}\subset\Siyz,\\
&\\
&{\rm deg}\,\sigma_J\,g_J,\:{\rm deg}\,\psi_K\,h_K\,\leq 2r\quad\forall\,J\subseteq\I_{xy},\,K\subseteq\I_{yz}
\end{array}\right..
\end{equation}
Equivalently, $\Q^*_r$ reads
\begin{equation}
\label{newformdual}
\Q^*_r:\quad
\displaystyle\max_{\lambda}\:\{\:\lambda \::\: f-\lambda\,\in\,P_r(g)+P_r(h)\:\},
\end{equation}
where $P_r(g)\subset P(g)$ denote the set of elements
$\sigma\in P(g)$ that can be written
as in (\ref{preordering}) and where in addition, 
${\rm deg}\,(\sigma_J\,g_J)\,\leq 2r$ for all $J\subseteq \I_{xy}$, and similarly for $P_r(h)$.

The first important step (i)  is to prove that $\inf\Q_r\uparrow f^*$ as $r\to\infty$, with
$f^*:=\min\:\{\,f(x,y,z)\,\vert\: (x,y,z)\in \K\,\}$, or equivalently,
\begin{equation}
\label{lower2}
f^*\,=\,\min _{\psi}\:\{\:\int f\,d\psi\:\vert\quad \psi(\R^{n+m+p}\setminus \K)=0;\:
\psi(\K)=1\},
\end{equation}
where the infimum is taken over all Borel probability measures $\psi$ 
on $\R^{n+m+p}$.

 The second important step (ii) is to prove absence of a duality gap between $\Q_r$ and 
its dual $\Q^*_r$. The final step (iii) easily follows from (a) and (b).\\

{\it Step (i)}.
We first prove that $\inf\Q_r\leq  f^*$ for all $r\geq r_0$.
Let $(x_0,y_0,z_0)\in \K$, and let $\mu:=\delta_{(x_0,y_0,z_0)}$ be the Dirac 
probability measure at $(x_0,y_0,z_0)$.
Let $u:=\{u_{\alpha\beta\gamma}\}$ be its (well defined) sequence of moments.
Then, obviously, $M_r(u,xy)\succeq0$, and $M_r(u,yz)\succeq0$.
Next, as $\mu$ is supported on $\K$, we obviously have
$M_r(g_J\,u,xy)\succeq0$, for every $J\subseteq \I_{xy}$ and $r\in\N$.
Similarly, $M_r(h_K\,u,yz)\succeq0$, for every $K\subseteq \I_{yz}$ and $r\in\N$.
Therefore, $u$ is feasible for $\Q_r$ with value $L_y(f)=\int fd\mu =f(x_0,y_0,z_0)$, which
proves that
\begin{equation}
\label{lower}
\inf\Q_r\,\leq\,f^*\,:=\,\min\:\{\:f(x,y,z)\:\vert\quad (x,y,z)\,\in\, \K\:\},\qquad\forall r\geq r_0.
\end{equation}

We next prove that $\inf\Q_r>-\infty$.
Let $r\geq r_0$ be fixed. As $\K_{xy}$ is compact, there is some $N_r$ such that the polynomials
$N_r\pm X^\alpha Y^\beta$, $\vert\alpha+\beta\vert\leq 2r$,
are all positive on  $\K_{xy}$.
By Schm\"udgen's Positivstellensatz \cite{schmudgen}, 
they all belong to the preordering $P(g)$. There is even some $l(r)$ such that
they all belong to $P_{l(r)}(g)$.
Similarly, as $\K_{yz}$ is compact, there is some $t(r)$ such that
the polynomials $N'_r-Y^\alpha Z^\beta$ belong to $P_{t(r)}(h)$, for 
some $N'_r$ and all $\vert\alpha+\beta\vert\leq 2r$. 

So, let $u^s$ be an arbitrary feasible solution of $\Q_s$, with
$s\geq \max[l(r),t(r)]$. One has
\[L_{u^s}(N_r\pm X^\alpha Y^\beta)\,=\,N_r\,u_0^s\,\pm\, u^s_{\alpha\beta\,0}\geq 0,
\quad\forall\,\vert\alpha+\beta\vert\leq 2r,\]
because $N_r\pm X^\alpha Y^\beta=\sum_{J\subseteq \I_{xy}}\sigma_Jg_J$
with ${\rm deg}\,\sigma_J\,g_J\leq 2s$,  and 
\[L_{u^s}(N_r\pm X^\alpha Y^\beta)\,=\,\langle \mathbf{\sigma_J},M_{r-r(g_J)}(g_J\,u^s,xy)\mathbf{\sigma_J}\rangle\geq 0,\]
where the latter inequality follows from $M_{s-r(g_J)}(g_J\,u^s,xy)\succeq0$; see (\ref{mom2}).

Similarly, one has
\[L_{u^s}(N'_r\pm Y^\alpha Z^\beta)\,=\,N'_r\,u^s_0\,\pm\, u^s_{0\,\alpha\beta}\geq 0,
\quad\forall\,\vert\alpha+\beta\vert\leq 2r,\]
because $N'_r\pm Y^\alpha Z^\beta=\sum_{K\subseteq \I_{yz}}\sigma_K\,h_K$
with ${\rm deg}\,\sigma_K\,h_K\leq 2s$,
and $M_{s-r(h_K)}(h_K\,u^s,yz)\succeq0$.
Therefore, 
\begin{equation}
\label{oui1}
\vert u^s_{\alpha\beta\,0}\vert \,\leq\,N_r\quad\forall \vert\alpha+\beta\vert\leq 2r;
\quad \vert u^s_{0\,\alpha\beta}\vert \,\leq\,N'_r\quad\forall \vert\alpha+\beta\vert\leq 2r.
\end{equation}
In particular,
\[\inf\Q_s\geq -\left[N_r\sum_{\alpha\beta}\vert (f_{xy})_{\alpha\beta\,0}\vert
+N'_r\sum_{\alpha\beta}\vert (f_{yz})_{0\,\alpha\beta}\vert
\right]\,>\,-\infty,\]
provided $s$ is sufficiently large.

So, for sufficiently large $s$, let
$u^s$ be a nearly optimal solution of $\Q_s$, i.e.,
\begin{equation}
\label{nearly}
\inf\Q_s\,\leq\,L_{u^s}(f)\,\leq\,\inf\Q_s+\frac{1}{s}
\,\leq\,f^*+\frac{1}{s},
\end{equation}
and complete each sequence $u^s$ with zeros to make it an
infinite sequence indexed in the 
canonical basis of $\R[X,Y,Z]$. Notice that by doing so,
only elements of the form $u^s_{\alpha\beta\,0}$ and $u^s_{0\,\alpha\beta}$ are non zero.
As (\ref{oui1}) is true for arbitrary $r$, by a standard diagonal argument, there exists a subsequence $\{s_k\}$ and an infinite sequence $u$ indexed in the 
canonical basis of $\R[X,Y,Z]$, such that
\begin{equation}
\label{convergence1}
u^{s_k}_{\alpha\beta\gamma}\,\to\, u_{\alpha\beta\gamma}\,\qquad\forall 
\alpha\in\N^{n},\quad\beta\in\N^{m},\quad\gamma\in\N^{p}.
\end{equation}
From what precedes the only non zero elements $u_{\alpha\beta\gamma}$ 
of $u$, are those with $\alpha=0$ or $\gamma =0$. Next, introduce the subsequences
\begin{equation}
\label{uk}
u^{xy}:=\,\{u_{\alpha\beta\gamma}\::\:\gamma\,=\,0\:\}\,;
\quad u^{yz}:=\,\{u_{\alpha\beta\gamma}\::\:\alpha\,=\,0\:\}.
\end{equation}
Recall that the matrix $M_r(u,xy)$ (resp. $M_r(u,yz)$) is identical to the moment matrix 
$M_r(u^{xy})$ (resp. $M_r(u^{yz})$) of 
the sequence $u^{xy}$ (resp. $u^{yz}$) indexed in the canonical basis of $\R[X,Y]$
(resp. $\R[Y,Z]$).

Similarly, as $g_J\in\R[X,Y]$ and
$h_K\in\R[Y,Z]$, the matrix $M_r(g_J\,u,xy)$ (resp. $M_r(h_K\,u,yz)$) is identical to the localizing matrix 
$M_r(g_J\,u^{xy})$ (resp. $M_r(h_K\,u^{yz})$) of 
the sequence $u^{xy}$ (resp. $u^{yz}$) indexed in the canonical basis of $\R[X,Y]$
(resp. $\R[Y,Z]$).

Next, let $r$ be fixed arbitrary.  Then from the feasibility of
$u^{s_k}$ in $\Q_{s_k}$, and the convergence
(\ref{convergence1}), we obtain
$M_r(u,xy)\,=\,M_r(u^{xy})\succeq0$ and $M_r(u,yz)=M_r(u^{yz})\succeq0$.
With same arguments, we also have
\begin{eqnarray*}
M_r(g_J\,u,xy)&=&M_r(g_J\,u^{xy})\succeq0,\quad\forall J\subseteq\I_{xy}\\
M_r(h_K\,u,yz)&=&M_r(h_K\,u^{yz})\succeq0,\quad\forall K\subseteq\I_{yz}.
\end{eqnarray*}

As $r$ was arbitrary, by Schm\"udgen's Positivestellensatz \cite{schmudgen}, it follows that the sequence $u^{xy}$ (resp. $u^{yz}$) has a representing measure
$\mu^{xy}$ (resp. $\mu^{yz}$)  with support contained in the compact set
$\K_{xy}$ (resp. $\K_{yz}$). Observe that
\[u^{xy}_{0\,\beta}\,=\,u_{0\,\beta \,0}\,=\,u^{yz}_{\beta\,0},\qquad\forall\,\beta\in\N^m.\]
Therefore, as measures on compact sets are moment determinate, $\mu^{xy}$ and $\mu^{yz}$ have same marginal $\mu^y$ on $\R^m$. 

Next, the probability measure $\mu^{xy}$ on the cartesian product of Borel spaces 
$\R^n\times\R^m$ can be disintegrated into
a stochastic kernel $q(dX\,\vert\, Y)$ on $\R^{n}$ given
$\R^m$, and its marginal $\mu^{y}$ on $\R^{m}$, i.e.,
\[\mu^{xy}(A\times B)\,=\,\int_B \,q(A\,\vert\, Y)\,\mu^y(dY),\]
for all Borel rectangles $(A\times B)$ of $\R^{n}\times\R^{m}$; see e.g.
Bertsekas and Schreve \cite[p. 139-141]{bertsekas}.
Similarly, the probability measure $\mu^{yz}$ on the cartesian product of Borel spaces 
$\R^m\times\R^p$
can be disintegrated into
a stochastic kernel $q'(dZ\,\vert\, Y)$ on $\R^{p}$ given
$\R^m$, and its marginal $\mu^{y}$ on $\R^{m}$, i.e.,
\[\mu^{yz}(B\times A)\,=\,\int_B \,q'(A\,\vert\, Y)\,\mu^y(dY),\]
for all Borel rectangles $(B\times A)$ of $\R^{m}\times\R^{p}$. 

Let $\mu$ be the probability measure on $\R^{n+m+p}$ defined by:
\begin{equation}
\label{defmu}
\mu(A\times B\times C)\,=\,\int_B \,q(A\,\vert\, Y)\,q'(C\,\vert\, Y)\,
\mu^y(dY),
\end{equation}
for all Borel rectangles $(A\times B\times C)$ of $\R^n\times \R^{m}\times\R^{p}$. 

Taking $C:=\R^p$ in (\ref{defmu}) yields that
\[\mu(A\times B\times \R^p)\,=\,
\int_B \,q(A\,\vert\, Y)\,\mu^y(dY)\,=\,\mu^{xy}(A\times B),\]
i.e., $\mu^{xy}$ is the marginal of $\mu$ on $\R^{n+m}$
(and in fact on $\K_{xy}$ because $\mu^{xy}(\K_{xy})=1$). 

Similarly, taking now $A=\R^n$, $\mu^{yz}$ is the marginal of $\mu$ on $\R^{m+p}$
(and in fact on $\K_{yz}$). This clearly implies that $\mu$ is supported on $\K$, i.e., $\mu(\K)=1$. Indeed,
\[\K\,=\,\{\:(x,y,z)\,:\:(x,y)\,\in\,\K_{xy}\,;\quad (y,z)\,\in\,\K_{yz}\:\},\]
and so,
\[\mu(\K)\,=\,\int q(\{x\::\,(x,y)\,\in\,\K_{xy}\}\,\vert\,Y)\,
q'(\{z\,:\,(y,z)\,\in\,\K_{yz}\}\,\vert\,Y)\,\mu^y(dY)\,=\,1,\]
because from the definitions of $\mu^y$, $\K_{xy}$ and $\K_{yz}$,
\[q(\{x\,:\,(x,y)\,\in\,\K_{xy}\}\,\vert\,Y)\,=\,
q'(\{z\,:\,(y,z)\,\in\,\K_{yz}\}\,\vert\,Y)\,=\,1\qquad\mu^y{\rm -a.e.}\]
Finally, observe that from the convergence (\ref{convergence1}), we obtain
\begin{equation}
\label{ouf}
\lim_{k\to\infty}\,L_{u^{s_k}}(f)\,=\,L_{u}(f)\,=\,\int f\,d\mu,
\end{equation}
because $f\in\R[X,Y]+\R[Y,Z]$. Therefore, by (\ref{nearly}) and
(\ref{ouf}), we get
\[f^*\,\geq\,\lim_{k\to\infty}\,\inf\Q_{s_k}\,=\,
\lim_{k\to\infty}\,L_{u^{s_k}}(f)\,=\,L_{u}(f)\,=\,\int f\,d\mu.\]
In view of (\ref{lower2}) and as $\mu$ is supported on $\K$, it follows that 
$f^*=\int f\,d\mu$. Therefore, $\inf\Q_{s_k}\to f^*$, and as the sequence $\{\Q_r\}$ is monotone nondecreasing, we obtain $\inf\Q_r\uparrow f^*$.

{\it Step (ii).} To prove absence of a duality gap between 
$\Q_r$ and its dual $\Q^*_r$, let $\nu$ be 
the uniform probability measure on $\K$, and let
$u$ be its sequence of moments.
As $\K$ has nonempty interior, it follows
that $M_r(g_J\,u)\succ0$, for every $J\subseteq \I_{xy}$, and all $r$, and similarly,
$M_r(h_K\,u)\succ0$, for every $K\subseteq \I_{yz}$, and all $r$.

But this implies that $M_r(g_J\,u,xy)\succ0$ for all $J\subseteq \I_{xy}$ and all $r$,
because $M_r(g_J\,u,xy)$ is a submatrix of $M_r(g_J\,u)$.
Similarly, $M_r(h_K\,u,yz)\succ0$ for all $K\subseteq \I_{yz}$ and all $r$,
because $M_r(h_K\,u,yz)$ is a submatrix of $M_r(h_K\,u)$.

This means that $u$ is a {\it strictly feasible} solution for $\Q_r$ and so,
Slater's condition is satisfied for $\Q_r$ (see e.g. \cite{boyd}). This in turn implies that there is
no duality gap between $\Q_r$ and its dual $\Q^*_r$, i.e.
$\sup\Q^*_r=\inf\Q_r$ for all $r\geq0$, and $\Q^*_r$ is even solvable 
(i.e. $\sup\Q^*_r=\max\Q^*_r$) if
$\Q_r$ has finite value; for more details on duality for semidefinite programs, see e.g. Vandenberghe and Boyd \cite{boyd}.

{\it Step (iii).} So let $f$ be strictly positive on $\K$, and let $f^*>0$ be its global minimum on $\K$, i.e., $f\geq f^*$ on $\K$. From (i)-(ii), there exists some $r$ such that
$\max\Q^*_r=\inf\Q_r\geq f^*/2>0$. Therefore, 
let $\lambda$ be an optimal solution of $\Q^*_r$. We have
$f-\lambda\in P(g)+P(h)$.  But then 
$f\in P(g)+P(h)$ because $\lambda>0$, the desired result.

(b) The proof of (b) is the same as that of (a), except that we now invoke 
Putinar's Positivstellensatz rather than Schm\"udgen's. Indeed,
if $N-\Vert (X,Y)\Vert^2\in Q(g)$, then 
Putinar's Positivstellensatz \cite{putinar} holds, i.e.,
every polynomial of $\R[X,Y]$, (strictly) positive on $\K_{xy}$, 
belongs to $Q(g)$; see also Jacobi and Prestel \cite{jacobi}. And so,
the polynomial $N_r\pm X^\alpha Y^\beta$, strictly positive on $\K_{xy}$
for sufficiently large $N_r$, belongs to the quadratic module $Q(g)$
(instead of the preordering $P(g)$). The rest of the proof is identical.

(c) The proof of (c) resembles that of (a), but the optimization problem
$\Q_r$ is now the {\it linear programming} problem
\begin{equation}
\label{pbplinear}
\L_r:\quad\left\{\begin{array}{ll}
\displaystyle\min_{u} &L_u(f)\\
\mbox{s.t.}&L_u(g^\alpha(1-g)^\beta)\geq0,\quad \forall \,(\alpha,\beta)\in G_r\\
&L_u(h^\alpha(1-h)^\beta)\geq0,\quad \forall\, (\alpha,\beta)\in H_r\\
&u_0=1.\\
\end{array}\right.,
\end{equation}
where:
\begin{eqnarray*}
G_r&=&\{(\alpha,\beta)\in\N^{\vert\I_{xy}\vert}\times\N^{\vert\I_{xy}\vert} \::\:
\sum_{j\in \I_{xy}}{\rm deg}\,g_j^{\alpha_j}+{\rm deg}\,g_j^{\beta_j}\,\leq\,2r\} \\
H_r&=&\{(\alpha,\beta)\in\N^{\vert\I_{yz}\vert}\times\N^{\vert\I_{yz}\vert} \::\:
\sum_{k\in \I_{yz}}{\rm deg}\,h_k^{\alpha_k}+{\rm deg}\,h_k^{\beta_k}\,\leq\,2r\}.
\end{eqnarray*}
The dual of the linear program $\L_r$ is the linear program
\begin{equation}
\label{newformduallp}
\L^*_r:\quad
\displaystyle\max_{\lambda}\:\{\:\lambda \::\: f-\lambda\,\in\,C_r(g,1-g)+C_r(h,1-h)\:\},
\end{equation}
where $C_r(g,1-g)\subset C(g,1-g)$ (resp. $C_r(h,1-h)\subset C(h,1-h)$)
is the subcone
\begin{eqnarray*}
C_r(g,1-g)&=&\{\,\sigma\in\R[X,Y]\::\:\sigma\,=\,\sum_{(\alpha,\beta)\in G_r}\,c_{\alpha\beta}\,g^\alpha\,(1-g)^\beta;\quad c_{\alpha\beta}\geq0\,\}\\
C_r(h,1-h))&=&\{\,\sigma\in\R[Y,Z]\::\:\sigma\,=\,\sum_{(\alpha,\beta)\in H_r}\,c_{\alpha\beta}\, h^\alpha\,(1-h)^\beta;\quad c_{\alpha\beta}\geq0\,\}.
\end{eqnarray*}
By compactness of $\K_{xy}$ and $\K_{yz}$, the polynomial $N_r\pm X^\alpha Y^\beta$ is strictly positive on $\K_{xy}$, for some $N_r$, and all $\alpha,\beta$ with $\vert\alpha+\beta\vert\leq 2r$. 
Therefore, using Krivine \cite{krivine} and Vasilescu \cite{vasilescu} instead of Schm\"udgen Positivstellensatz, there is some $l(r)$ such that 
 $N_r\pm X^\alpha Y^\beta\in C_{l(r)}(g,1-g)$.
Proceeding as in (a), the latter property is used in an optimal solution $u^s$ 
of $\L_s$, to bound $\vert u^s_{\alpha\beta 0}\vert$ for all $\vert\alpha+\beta\vert\leq 2r$, uniformly in $s$.
We then obtain the convergence
(\ref{convergence1}) for a subsequence $\{s_k\}$.
Using again Krivine \cite{krivine} and Vasilescu \cite{vasilescu},
the subsequence $u^{xy}$ of the limit sequence $u$ (see (\ref{uk})), is the moment vector of a measure $\mu^{xy}$ supported on $\K_{xy}$, and 
$u^{yz}$ is the moment vector of a measure
$\mu^{yz}$ supported on $\K_{yz}$. In the present case, 
from linear programming duality,
there is no duality gap between the linear programs $\L_r$ and $\L^*_r$. 
The rest of the proof is along the same lines.

\subsection{Proof of Theorem \ref{coro1}}
\label{proofcoro1}
As $\K_x$ and $\K_{yz}$ are both compact, then (possibly after 
some change of variables), we may and will assume that $\K_x\subset [-1/2,1/2]^n$ and $\K_{yz}\subset [-1/2,1/2]^{m+p}$.
We only prove (a) because similar arguments hold for (b).
Introduce the optimization problem:
\begin{equation}
\label{pbpnew}
\Q_r:\quad\left\{\begin{array}{ll}
\displaystyle\min_{u} &L_u(f)\\
\mbox{s.t.}&M_r(u,xy)\succeq0\\
&M_{r-r(g_J)}(g_J\,u,x)\succeq0,\quad \forall\, J\subseteq\I_{x}\\
&M_{r-r(h_k)}(h_K\,u,yz)\succeq0,\quad \forall\, K\subseteq\I_{yz}\\
&u_0=1.\\
\end{array}\right.,
\end{equation}
where now $M_r(g_J\,u,x)$ is the obvious analogue of
$M_r(g_J\,u,xy)$ defined in \S \ref{notation}. The dual $\Q^*_r$ now reads
\begin{equation}
\label{pbpdualnew}
\Q^*_r:\quad\left\{\begin{array}{ll}
\displaystyle\max_{\{\sigma_J\},\{\psi_K\},\kappa_\emptyset,\lambda} &\lambda\\
\mbox{s.t.}&
f-\lambda\,=\,\kappa_\emptyset+\displaystyle\sum_{J\subseteq\I_{x}}\,\sigma_J\,g_J
+\displaystyle\sum_{K\subseteq\I_{yz}}\,\psi_K\,h_K\\
&\\
&\kappa_\emptyset\in\Sixy,\:\{\sigma_J\}\subset\Six
,\quad\{\psi_K\}\subset\Siyz,\\
&\\
&{\rm deg}\,\kappa_\emptyset,\:{\rm deg}\,\sigma_J\,g_J,\:{\rm deg}\,\psi_K\,h_K\,\leq 2r
\end{array}\right..
\end{equation}
Equivalently, $\Q^*_r$ reads
\begin{equation}
\label{newformdualnew}
\Q^*_r:\quad
\displaystyle\max_{\lambda}\:\{\:\lambda \::\: f-\lambda\,\in\,(\Sixy)_r+P_r(g)+P_r(h)\:\},
\end{equation}
where $(\Sixy)_r$  is the set of elements of $\Sixy$ of degree at most $2r$.

Proceed as in the proof of Theorem \ref{thmain}(a), and so, let
$u^s=(u^s_{\alpha\beta\gamma})$ be an arbitrary feasible solution of $\Q_s$ in (\ref{pbpnew}).
As $\K_{yz}\subset [-1/2,1/2]^{m+p}$, $1\pm Y^\alpha Z^\beta \in P(h)$ for all
$(\alpha,\beta)\in\N^m\times\N^p$. Similarly, As $\K_{x}\subset [-1/2,1/2]^{n}$,
$1\pm X^\alpha\in P(g)$ for all $\alpha\in \N^n$. Therefore, let $r$ be fixed, arbitrary. With same arguments as in the proof of Theorem \ref{thmain}(a), for all sufficiently large $s$,
\begin{eqnarray}
\label{ehoui1}
\vert u^s_{\alpha 00}\vert <1&&\forall\,\alpha\in\N^n;\quad\vert\alpha\vert\leq 2r\\
\label{ehoui2}
\vert u^s_{0\beta\gamma}\vert <1&&\forall\,(\beta,\gamma)\in\N^m\times\N^p,\quad\vert\beta+\gamma\vert\leq 2r.
\end{eqnarray}
There is some additional technicality because we also need the boundedness of $\vert u^s_{\alpha\beta0}\vert$, uniformly in $s$. But 
from $M_{s}(u,xy)\succeq0$, we get $u^{s}_{(2\alpha) 00}\,u^{s}_{0(2\beta)0}\geq (u^{s}_{\alpha\beta0})^2$, and so, in view of (\ref{ehoui1})-(\ref{ehoui2}),
$1>\vert u^s_{\alpha\beta0}\vert$ for all $s$ sufficiently large. 
Therefore, as in the proof of Theorem \ref{thmain}(a) let $\{u^s\}$
be a sequence of nearly optimal solutions of $\Q_s$.
The limit sequences $u^{xy}$ and $u^{yz}$ in (\ref{uk}) 
satisfy $\vert u^{xy}_{\alpha\beta}\vert\leq 1$
and $\vert u^{yz}_{\beta\gamma}\vert\leq 1$ for all $\alpha,\beta,\gamma$.
But this implies that $u^{xy}$ is the moment sequence of a probability
measure $\mu^{xy}$ supported on $[-1,1]^{n+m}$; see Berg \cite[Theor. 9]{berg}. 
In addition, again as in the proof of
Theorem \ref{thmain}(a), $u^{x}$ and $u^{yz}$ are moment sequences of two probability measures
$\mu^x$ and $\mu^{yz}$ with support contained in $\K_x$ and $\K_{yz}$ respectively.
By construction,
\begin{equation}
\label{compact}
\int X^\alpha\,d\mu^x\,=\,u^x_{\alpha00}\,=\,
u^{xy}_{\alpha 00}\,=\,\int X^\alpha\,d\mu^{xy},\qquad\forall\,\alpha\in\N^n.
\end{equation}
As both $\mu^{xy}$ and $\mu^x$ have their support in a compact set, they are moment determinate, and so (\ref{compact}) implies that $\mu^x$ is the marginal of $\mu^{xy}$ on $\K_x$. Similarly, the marginal $\mu^y$ of $\mu^{yz}$ on $\R^m$ is the same 
as the marginal of $\mu^{xy}$ on $\R^m$. Therefore, the measure
$\mu$ defined in (\ref{defmu}) has marginal $\mu^x$ on $\K_x$, marginal
$\mu^{yz}$ on $\K_{yz}$, and marginal $\mu^{xy}$ on $\R^n\times\R^m$.
The rest of the proof is the same as that of Theorem \ref{thmain}(a).

(b) One proceeds exactly as in (a), except that now one invokes
Putinar's instead of Schm\"udgen's Positivstellensatz.

\section*{Acknowledgements}
This work was done under (french) ANR grant NT05-3-41612.
The author wishes to thank M. Schweighofer and T. Netzer for helpful remarks 
and suggestions for \cite{largescale}, that we have also used in the present paper.

\end{document}